\documentclass[12pt,english]{article}
\usepackage{geometry}
\usepackage{color}
\usepackage{amsmath}
\usepackage{amssymb}
\usepackage{amsthm}

\def\E{\hskip.15ex\mathrm{E}\hskip.10ex}
\def\P{\mathrm{P}}

\def\phi{\varphi}



\usepackage{babel}

\newtheorem{theorem}{Theorem}
\newtheorem{lemma}[theorem]{Lemma}

\newtheorem{corollary}[theorem]{Corollary}

\newtheorem{example}[theorem]{Example}

\newif\ifproofs
\proofstrue

\newif\ifp
\pfalse

\begin{document}
\global\long\def\E{\mathbb{E}}
\global\long\def\P{\mathbb{P}}
\global\long\def\N{\mathbb{N}}
\global\long\def\ind{\mathbb{I}}

\title{
{\normalsize\tt\hfill\jobname.tex}\\
On Markov--up processes and their recurrence properties 
}
\author{A.Yu. Veretennikov\footnote{Institute for Information Transmission Problems;
email: ayv @ iitp.ru. Supported by 
grant 
Russian Foundation of Basic Research 20-01-00575a 
(Theorem \ref{thm1}, Lemma \ref{lem3})}, 
M.A. Veretennikova\footnote{University of Warwick, Coventry, United Kingdom;  email: maveretenn @ gmail.com. 
}
}

\maketitle
\begin{abstract}
A simple model of the new notion of ``Markov up'' processes is proposed; its positive recurrence and ergodic properties are shown under the appropriate conditions. 

\medskip


\medskip


\end{abstract}


\section{Introduction: the model}\label{sec:intro}
The idea of integer valued processes which behave like markovian on the periods of growing and in a more complicated non-markovian way on the periods of decreasing was suggested by Alexander Dmitrievich Solovyev in a private communication in the late 90's \cite{ADSprivate}. In this period of his research activity he only worked on applied projects; hence, there is no doubt that this idea was also an applied one, most likely related to the theory of reliability.  To the best of the authors knowledge he did not leave any notes on this theme. Also, the authors are not aware of any publications on this topic, although certain close models do exist in the literature. In this paper  a toy model of this idea is proposed. 

Consider a process $X_n, n\ge 0$  on  $\mathbb Z_+ = \{0, 1, \ldots\}$, or on  $\mathbb Z_{0,\bar N}=\{0, \ldots \bar N\}$ with some $0<\bar N < \infty$, possessing the following property: for any  $n\ge 1$ where the last jump was up (including staying), it is assumed that for some function $\phi(i,j), i,j\in \mathbb Z_+$, 
\begin{equation}\label{e1}
\mathbb P(X_{n+1} = j | {\cal F}^X_n; X_n\ge X_{n-1}) = 
\phi(X_n,j), 
\end{equation}
that is, the ``movement upwards remains markovian''; the ``decision'' to turn downwards is also markovian in the first instant; 
however, where the last jump was down, the next  probability distribution may depend on some part of the past trajectory: namely for some function $\psi(X_n, \ldots X_{\zeta_n})$
\begin{equation}\label{e2}
\mathbb P(X_{n+1} = j | {\cal F}^X_n; X_n <  X_{n-1}) = 
\psi(X_n, \ldots X_{\zeta_n}),
\end{equation}
where $\zeta_n$ is the last turning time from ``up'' to ``down'' before $n$; it is formally defined  in (\ref{zeta}) in what follows. 
Hence, the ``memory'' of the process while moving down 
is limited by the last time of turning down; the latter moment may not be bounded.
These assumptions  reflect the property that while the process ``goes up" its transition probabilities for any jump up obey the Markov property (\ref{e1}); as soon as it goes down, its transition probabilities ``must" remember some past values of the trajectory, namely, from the last jump up moment. The case of equality $X_n = X_{n-1}$ -- or staying at its place -- is included in the movement up; probably it {\em may} be shifted to the movement down, but apparently it would change the calculus and we do not pursue to study all possibilities here at once. Also, some more complicated rules could be introduced instead of those described above; however, our goal is just to show a simplest version of the idea of a ``Markov up" process and to discuss some recurrence and ergodic properties which this model may possess. 

As a rationale, the model may be applied to a situation of the evolution of some involved many-component device  which may have several states and which ``goes up" while it is working , or it ``goes down" if one or several of the critical components in this device break down, after which the evolution does not stop but becomes more and more chaotic with a likely further disbalance or even a, or, at least, dependent on all the states after their break down: the device ``remembers" the event of the faults in the critical components all the time until they are repaired (in the simplest example fixing is just reloading the system), after which the transmission may resume again and the behaviour becomes ``markovian" again, satisfying the condition (\ref{e1}). There is also some evidence that certain disastrous processes related to complicated devices may expose similar features: once some critical failure occurs, the process of destruction may accelerate and be unpredictably chaotic until some rescue arrives.

Note that the probability of such a model to be at some subset in the state space could be viewed as a characteristic of the reliability of this device. Suppose that the movement ``up'' of the process $X$ is treated as approaching to some goal which is a high enough level above zero, and that at any moment of time the position of $X$ above the minimal level $N$ brings some profit, while falling down below the level $N$ is regarded as a failure with no dividends or even with some loss due to the expenses for repairing with the necessity to recover and to start raising up again. Then  the dynamical, or instantaneous reliability of the system may be defined as the probability $r(t):={\mathbb P}(X_t>N)$. Naturally, we are interested in computing this function $r(t)$, or, at least, its limit $r(\infty):=\lim_{t\to\infty}r(t)$, or its stationary value if the latter exists.  Indeed, traditionally all features of a model are evaluated and described in a stationary regime. It is well known that very often in probability models such a limit coincides with the stationary value of $r(t)$. In such a setting the property of a positive recurrence may help to show that this invariant or limiting probability $r(\infty)$ exists. The next important question would be to find the rate of this convergence; it is not pursued in this paper. The issue of the bounds for the rate of this convergence is left until further research and publications. Here we just recall that positive recurrence is naturally linked to the existence of a stationary regime (see the corollary \ref{cor} in what follows), and the stronger recurrence, the faster convergence.


The paper \cite{Martell} proposes a Markov model for the daily dynamics of the Fire Weather Index (FWI), which estimates the risk of wildfire. The authors do indicate that in fact the probability of wildfire escaping will grow as the duration of a several-day intensive fire onset increases. Statistical analysis in the paper concerns the suitable order of the Markov chain. It shows that for the analyzed data mostly a Markov chain of order 1 is suitable, however sometimes order 2 is preferable. Data is limited to the province of Ontario, and the appropriate order may be different elsewhere. In our model the length of memory is not fixed, which allows greater flexibility. It also takes into account the duration of last fire onset, which may be beneficial. For example, the paper \cite{Ramachandran} 
supports the idea that the total area burnt by a fire is an exponential function of time after ignition. Such amplification of chaos and further imbalance is discussed in the previous paragraph about the functionality of a multi-component device. Evidence of local memory dependence suggests that possibly a Markov-up process should be a reasonable model for evolution of an index which quantifies realistic damage from fire. For the process to be called Markov-up the worse the prognosis of the total damage the lower the index should be. In case of working with a variable such as FWI ranging from 0 (low danger) to 100 (extreme danger), perhaps, we could just as well introduce the notion of a ``Markov-down'' process, reversing the directions of jumps with the specified transition probability characteristics. For fire damage index dynamics it would be appropriate to consider a variation of the Markov-down process, in which return to Markov behaviour happens after the index reaches a `low' danger threshold level in several sequential steps.  Note that this index may also be regarded as a reliability type characteristic where the reliability value could be defined as a probability that this index does not exceed some level. What is more, actually, the probability of each possible value of this index could be a more accurate and informative characteristic of an ``extended reliability'' type.  The theory in this paper concerns the simplest version of a Markov-up process.

\ifproofs
Note that according to  (\ref{e2}) the ``transition probabilities" $\mathbb P(X_{n+1} = j | X_n, \ldots  X_{\zeta_n})$ after jumps down do not depend on $n$, that is, 
\begin{align*}
\mathbb P(X_{n+1} = j | X_n, \ldots,  X_{\zeta_n})|_{\zeta_n=m, X_n=a_0, \ldots,  X_{\zeta_n}=a_m} 
 \\\\
= \mathbb P(X_{n+k+1} = j | X_{n+k}, \ldots,  X_{\zeta_{n+k}})|_{\zeta_{n+k}=m, X_{n+k}=a_0, \ldots,  X_{\zeta_{n+k}}=a_m},
\end{align*}
for any $m\ge 0$ and $k\ge 0$ in the case of 
$$
a_0 > \ldots >a_m, 
$$
where it is assumed that $X_{\zeta_n -1}\le a_m$.
The similar assumption is made about the probabilities $\mathbb P(X_{n+1} = j |X_n)$ after jumps up, see (\ref{e1}). This corresponds to the ``homogeneous" situation, in which it makes sense to pose a question about ergodic properties.   
For the conditional probabilities after the ``jumps down" the memory could be, in principle, unlimited, in the sense that it is not described by, say, m-Markov chains (i.e., with the memory of length $m$) except for the case of a finite $\bar N$. However, the process ``does not remember anything which is older than the last turn down", that is, there is no dependence of future probabilities on the past earlier than time $\zeta_n$ for each $n$. The moment $\zeta_n$ itself is interpreted as the last jump up before the fault occurs, and all the time before the faulty component is fixed, the device keeps record of what has happened from that moment to the present time, and the transition probabilities depend on this memory. The first jump up after a series of jumps down signifies that the faulty component is fixed and, hence, movement up resumes.  The movement in both directions can have several options, that is, it is not assumed that any jump up is by +1 and any jump down is with -1. Naturally, from zero there are only jumps up, or the process may stay at its place. The model with a finite $\bar N$ does not differ too much from the infinite version: since we are interested in bounds which would not depend on $\bar N$, the calculus would be very similar: the only point is that at $\bar N$ it should be specified what kind of jumps are possible; we do not pursue this version here assuming $\bar N = \infty$.

\fi

Models with more involved dependancies are possible: for example, instead of the immediate switching to ``Markov" probabilities after one jump up, it could be assumed that such a switch occurs after several steps up, or after the average in time of consequent jumps up or down exceeds some level, etc. Probably, some other adjustments of the model may be performed in order to include some specific forest fire features  mentioned earlier.


We are interested in establishing  ergodic properties for the model (\ref{e1})--(\ref{e2}) under certain ``recurrence" and ``non-singularity" assumptions. So, recurrence is one of the key points addressed here.

There are some ideological similarities of the proposed model with renewal processes, 
and with a (more general) notion of Hawkes processes, and also with semi-Markov processes.  Actually, this is a special case of semi-Markov type, as well as a special case of a regeneration process. Moreover, as we shall see in what follows, some transformation of the model based on the enlarged state space turns out to be a particular Markov process, which is not really surprising since, as is well-known, any process may be regarded as Markov after a certain change of the state space. Yet, this is not always useful. In any case, ergodic properties of the model are to be established from scratch, and markovian features will only be used in what concerns the invariant measure via the Harris -- Khasminskii principle. 

Because of many new objects, quite a few definitions will be repeatedly reminded to the reader during the text.

\section{Main results}

\noindent
We use standard notation $a\wedge b = \min(a,b)$, $a\vee b = \max(a,b)$.

~

\noindent
{\bf Further notations:}
Let us define for each $n \ge 0$ the random variables 
\begin{equation}\label{zeta}
\zeta_n:= \inf_{}(k \le n: \, \Delta X_i := X_{i+1} - X_{i}< 0, \forall i = k, \ldots n), \; (\inf(\emptyset) = +\infty),
\end{equation}
\begin{equation}\label{xi}
\xi_n : = \sup(k\ge n: \text{all increments}\; \Delta X_i\ge 0, \, \forall \, n\le i\le k)\vee n.
\end{equation}
\begin{equation}\label{chi}
\chi_n : = \sup(k\ge n: \text{all increments}\; \Delta X_i < 0, \, \forall \, n\le i\le k)\vee n.
\end{equation}
Also, let
\begin{align}\label{tf}
\hat X_{i,n}:= X_i 1(\zeta_n \wedge n\le i \le n), \quad 
\tilde {\cal F}_n 
= \sigma(\zeta_n; \hat X_{i,n}: 0\le i \le n).
\end{align}
Note that the family  $(\tilde {\cal F}_n)$ is not a filtration, and this is not required. We have, $\tilde {\cal F}_n \subset {\cal F}_n$ and $1(\zeta_n\wedge n = n)E(\xi|\tilde {\cal F}_n) = 1(\zeta_n\wedge n = n)E(\xi|X_n)$ $\forall \xi$. Also, note that $\hat X_{n,n} = X_n$ for any $n$.

~

\noindent
Now let us state the {\bf assumptions} which rewrite from scratch the  formulae (\ref{e1}) and (\ref{e2}).

~

\noindent
{\bf A1. Random memory depth:} 
{\em For any $n$, 
\begin{equation}\label{ctF}
P(X_{n+1}=j|{\cal F}_n) = P(X_{n+1}=j|\tilde {\cal F}_n) \quad \text{a.s.},
\end{equation}
and the latter conditional probability does not depend on $n$ given the past $X_n, \ldots X_{\zeta_n\wedge n}$, which serves as the analogue of the homogeneity.}\\
The random memory depth is what clearly distinguishes the proposed model from Markov chains with a fixed memory length also known as complex Markov chains.

\noindent
{\bf A2. Irreducibility (local mixing):} 
{\em For any $x\le N$ and for two states $y=x$ and $y=x+1$}
$$
P(X_{n+1} = y|\tilde {\cal F}_n, X_n=x) \ge \rho >0.
$$
%
Note that $2\rho \le 1$.
Along with the recurrence condition, the assumption A2 will guarantee the irreducibility of the process in the extended state space where the process becomes Markov, see  (\ref{markovn}) below.

\noindent
{\bf A3. Recurrence-1}: 
{\em There exists $N \ge 0$ such that 
\begin{equation}\label{q0}
P(\text{jump down}\equiv (X_{n+1}< X_n)|\tilde {\cal F}_n, N<X_n) \ge \kappa_0  >0;
\end{equation}
$$
P(
X_{n+1}< X_n
|\tilde {\cal F}_n, N< X_n < X_{n-1}) \ge \kappa_1>0, 
$$
etc., and for any $n\ge m$
\begin{equation}\label{q}
P(X_{n+1}< X_n
|\tilde{\cal F}_n, N<X_n < \ldots <X_{n-m+1}) \ge \kappa_{m-1}>0, 
\; \forall \, 1\le m,
\end{equation}
%
Note that $\kappa_0\le \kappa_1\le \ldots$}
Denote $$q=1-\kappa_0; \quad q<1.$$ Then
$$
P(\text{jump up}\equiv (X_{n+1}\ge X_n)|\tilde {\cal F}_n, N<X_n) \le 1- \kappa_0  = q < 1.
$$

~

\noindent
{\bf A4. Recurrence-2}:
{\em It is assumed that  the following infinite product converges
\begin{align}\label{prod}
\bar \kappa_\infty:=\prod_{i=0}^\infty \kappa_i >0; 
\end{align}
and 
\begin{align}\label{ka}
\sum_{i\ge 1} i (1-\kappa_i) <\infty.
\end{align}
Let 
$$
\bar q := 1-\bar \kappa_\infty (< 1) \quad \& \quad 
q: = 1 - \kappa_0 \quad (<1).
$$
Note that 
$$
P(\mbox{\em jump up}\equiv (X_{n+1}\ge  X_n)|\tilde {\cal F}_n, N<X_n) \le 1-  \kappa_0  = \bar \kappa_0 = q <1.
$$
}





~

\noindent
{\bf A5. Jump up moment bound:} 
\begin{align}\label{2moment}
M_1:= \mathop{\mbox{ess}}\limits_{P}\sup\limits_{
\omega} \sup\limits_{n} E ((X_{n+1}-X_{n})_+|\tilde {\cal F}_n) <\infty. 
\end{align}
Let $$\bar q= 1-\bar \kappa_\infty \; (<1).$$ 
This is the upper bound for the probability that the fall down is not successful, i.e., that the ``floor'' $[0,N]$ is not reached in one go.

~

Denote 
$$
\bar \kappa_m:=\prod_{i=0}^m \kappa_i \quad (\ge \bar \kappa_\infty >0). 
$$
Let us emphasize that the index $i$ in $\kappa_i$ is not the state where the process $X$ is, but the value for how long the process is falling down. 
The process remembers for how long it has been going down so far, and the longer it goes down the more probable is to continue in this direction, at least, until the process reaches $[0,N]$. Equivalently, 
$$
\sum \ln \kappa_i <\infty.
$$
Of course, this implies that $\kappa_i \to 1$ as $i\to\infty$, which is, clearly, a weaker condition than (\ref{prod}). Convergence of the sequence $\kappa_i$ to $1$, if it is monotonic, may be interpreted in a way that the longer is the decreasing trajectory, the more faulty components in the device: each jump down makes some additional disorder in the system, which further increases the probability to continue falling down.\\

~

\begin{example}
The assumption (\ref{ka}) is satisfied, for example, under the condition $1-\kappa_m \le \frac{C}{m^3}$, or, equivalently, 
$$
\kappa_m \ge 1-\frac{C}{m^3}.
$$
An exponential rate of the approach of the sequence $\kappa_m$ to $1$ accepted in some applied models of a fire evolution could be interpreted as the inequality
$$
\kappa_m \ge 1- \exp(-\lambda m)
$$
with some $\lambda>0$.

~

The assumption (\ref{2moment}) is valid, for example, if there exists a nonrandom constant $C\ge 0$ such that with probability one
$$
X_{n+1}-X_{n} \le C<\infty.
$$

\end{example}

~

Denote 
$$
\tau=\tau^1:= \inf(t\ge 0: X_t\le N); \; \gamma:=\inf(t\ge \tau: X_{t-1}\le X_t = N).
$$
The 
regeneration occurs not at moment $\tau$, but at moment $\gamma$. However, the expectation of $\gamma$ may be evaluated via $E_x\tau$.  
Hence, it will be useful to introduce by induction  the following two sequences of stopping times with respect to the filtration ${\cal F}^X_n$: 
\begin{align*}
 T^{n}:= \inf(t> \tau^n: X_t > N), \quad 
 \tau^{n+1}:= \inf(t > T^n: X_t\le N).
\end{align*}

~

\noindent
{\bf The convention.} With the initial position $X_0 = x$ we assume that any artificial ``admissible past'' is allowed, that is, we accept that there is some fictitious past which could have preceded this state;  we include in this past nothing if the artificial state $X_{-1}$ does not exceed $X_0$, or we add the fictitious past trajectory from the last starting moment of the fall $\zeta_0$: $X_{\zeta_0}, \ldots, X_{-1}$. 
From the assumption (A1) it follows that the process $(Y_n, {\cal F}^Y_n)$ is Markov; of course,  ${\cal F}^Y_n = {\cal F}^X_n$.

~

\ifproofs

Let us recall the definitions of Greeks: 
$$
\zeta_n:= \inf_{}(k \le n: \, \Delta X_i := X_{i+1} - X_{i}< 0, \forall i = k, \ldots n) \qquad (\inf(\emptyset) = \infty);
$$
$$
\hat X_{i,n}:= X_i 1(\zeta_n \wedge n)\le i \le n), \quad 
\tilde {\cal F}_n 
= \sigma(\zeta_n; \hat X_{i,n}: 0\le i \le n);
$$
$$
\xi_n : = \sup(k\ge n: \text{all increments}\; \Delta X_i\ge 0, \, \forall \, n\le i\le k)\vee n;
$$
$$
\chi_n : = \sup(k\ge n: \text{all increments}\; \Delta X_i < 0, \, \forall \, n\le i\le k)\vee n.
$$

\section{Auxiliary lemmata}


\begin{lemma}\label{lem1}
Under the assumption (A3) for any $x>N$, 
$$
E_x (\xi_0 - 0) \le M_2 = \frac{q}{(1- q)^2}.
$$
\end{lemma}
\noindent
{\bf Proof.} 
Recall that the random variable $\xi_n$ was defined by the formula
$$
\xi_n : = \sup(k\ge n: \text{all increments}\; \Delta X_i\ge 0, \, \forall \, n\le i\le k)\vee n.
$$
We use the notations from the proof of lemma \ref{lem3} (below): for $i\ge n$ let
$$
e_i = 1(X_{i+1}\ge X_i), \; \bar e_i = 1(X_{i+1}<X_i), \; \Delta X_i = X_{i+1} - X_i, \; \ell_n^i = \bar e_i \prod_{k=n}^{i-1} e_k  \; (\text{assume}\; \prod_n^{n-1}=1).
$$ 
The bounds in this lemma and in the other lemmata will not depend on the initial state $x$, so we drop this index in $E_x$ and $P_x$ in this section (but not in the proof of the main result).  
We have, for $i\ge n$
\begin{align*}
E_x(e_i|X_i>N) = P(X_{i+1}\ge X_i|X_i>N)
 \\\\
=E_x (P_x(X_{i+1}\ge X_i|\tilde {\cal F}_i, X_i>N)|X_i>N)
\le 1-\kappa_0=q.
\end{align*}
Then almost surely
$$
\xi_n - n  = \sum_{k=0}^{\infty} k\bar e_{n+k}\prod_{i=n}^{n+k-1}e_i = \sum_{k=1}^{\infty} k\bar e_{n+k}\prod_{i=n}^{n+k-1}e_i =  \sum_{k=1}^{\infty} k\ell_n^{n+k}.
$$
So, we estimate, 
\begin{align*}
E_{x} (\xi_n - n)  = E_{} \sum_{k=1}^{\infty} k\bar e_{n+k}\prod_{i=n}^{n+k-1}e_i \le \sum_{k=1}^{\infty} k E_{x}\prod_{i=n}^{n+k-1}e_i 
 \\
\le \sum_{k=1}^{\infty} k  q^{k} 
= q \sum_{k=1}^{\infty} k  q^{k-1} 
= \frac{q}{(1-q)^2}:=M_2 < \infty. \qquad \hfill \text{QED}
\end{align*}

~


\noindent
Let us recall, 
$$
\tau:= \inf(t\ge 0: X_t\le N), \quad \chi_n : = \sup(k\ge n: \text{all increments}\; \Delta X_i < 0, \, \forall \, n\le i\le k)\vee n,
$$
and
$$
P_x(X_{n+1}< X_n
|\tilde{\cal F}_n, N<X_n < \ldots <X_{n-m+1}) \ge \kappa_{m-1}>0, 
\quad \forall \, 1\le m,
$$
and also
$$
\chi_n : = \sup(k\ge n: \text{all increments}\; \Delta X_i < 0, \, \forall \, n\le i\le k)\vee n.
$$
\begin{lemma}\label{lem2}
Under the assumptions (A3)-(A4), for any $x>N$, ($n=0$)
$$
E_{x} (\chi_n - n)1(\chi_n<\tau) \le \sum_{i\ge 1} i (1-\kappa_{i})
: = M_3 < \infty.
$$
\end{lemma}
\noindent
{\bf Proof.} 
Similarly to the calculus of the previous lemma but with the replacement of $e_i$ by $\bar e_i$ and vice versa, we have
\begin{align*}
(\chi_n - n)1(\chi_n<\tau) \le  
\sum_{k=1}^{\infty} k e_{n+k} 1(n+k-1 <\tau)\prod_{i=n}^{n+k-1} \bar e_i, 
\end{align*}
so, 
\begin{align*}
E_x (\chi_n - n)1(\chi_n<\tau) \le  
E_x \sum_{k=1}^{\infty} k e_{n+k} 1(n+k-1 <\tau) \prod_{i=n}^{n+k-1} \bar e_i 
 \\
\le \sum_{k=1}^{\infty} k E_x 1(n+k-1 <\tau) (\prod_{i=n}^{n+k-1}\bar e_i) 
E_x(e_{n+k} | \Delta X_i<0, 0\le i\le n+k-1) 
 \\
\stackrel{A3}\le \sum_{k=1}^{\infty} k E_x 1(n+k-1 <\tau) (\prod_{i=n}^{n+k-1}\bar e_i) (1-\kappa_{k}) 
\le \sum_{k=1}^{\infty} k (1-\kappa_{k})
=: M_3 \stackrel{A4}<\infty. \quad \text{QED}
\end{align*}

~

\noindent
Let us recall once more, 
$$
\xi_n : = \sup(k\ge n: \text{all increments}\; \Delta X_i\ge 0, \, \forall \, n\le i\le k)\vee n.
$$

\begin{lemma}\label{lem3}
Under the assumptions (A3) and (A5) the expected value of the maximum positive increment over any single period of running up (non-strictly) until the first jump down is finite:
\begin{align*}
\sup_{n,x} E_x (X_{\xi_n} - X_{n})_+ 
 \le M_4 < \infty.
\end{align*}

\end{lemma}
\noindent
{\bf Proof.} 
\ifp
{\color{red}($p>1$).} Denote $e_i = 1(X_{i+1}\ge X_i)$, $\bar e_i = 1(X_{i+1}<X_i)$, $\Delta X_i = X_{i+1} - X_i$, $\ell_n^i = \bar e_i \prod_{k=n}^{i-1} e_k$. By H\"older's inequality
\begin{align*}
E(X_{\xi_n} - X_{n})^r_+ = E\sum_{i=n}^{\infty} \ell_n^i (X_i - X_n)^r
 \\
\le \sum_{i=n}^{\infty} (E \ell_n^i)^{1-r/p} (E\ell_n^i(X_i - X_n)^p)^{r/p}
\end{align*}
Further, 
\begin{align*}
E\ell_n^i(X_i - X_n)^p = E(\ell_n^i\sum_{j=n}^{i-1}\Delta X_j)^p
 \le E(\sum_{j=n}^{i-1} e_j\Delta X_j)^p 
 \\
= (i-n)^{p}E(\frac1{i-n}\sum_{j=n}^{i-1} e_j\Delta X_j)^p 
\stackrel{Jensen}{\le} (i-n)^{p}E\frac1{i-n}\sum_{j=n}^{i-1} (e_j\Delta X_j)^p 
 \\
= (i-n)^{p-1}\sum_{j=n}^{i-1} E(e_j\Delta X_j)^p 
\le (i-n)^{p-1}\sum_{j=n}^{i-1} M_1
= (i-n)^{p} M_1.
\end{align*}
Also, 
\begin{align*}
E \ell_n^i \le \bar\kappa_0^{i-n} \quad (<1) \quad (or \; \le \bar\kappa_0^{i-n-1}?)
\end{align*}
Hence, 
\begin{align*}
E(X_{\xi_n} - X_{n})_+ 
\le \sum_{i=n}^{\infty} (E \ell_n^i)^{1-r/p} (E\ell_n^i(X_i - X_n)^p)^{r/p} 
 \\
\le \sum_{i=n}^{\infty} (\bar\kappa_0^{i-n})^{1-r/p} ((i-n)^{p} M_1)^{r/p} =:M_4 < \infty,
\end{align*}
as required. Lemma \ref{lem3} for {\color{red}($p>1$)} is proved.

~

\fi
First of all, it suffices to show that
\begin{align*}
\sup_{n,x} E_x (X_{\xi_n} - X_{n})_+ |\tilde {\cal F}_n)
 \le M_4 < \infty.
\end{align*}
Further, we have
\begin{align*}
\sup_{n,x} E_x ((X_{\xi_n} - X_{n})_+ |\tilde {\cal F}_n, X_n\le N)
\le N+ \sup_{n,x} E_x ((X_{\xi_n} - X_{n})_+ |\tilde {\cal F}_n, X_n > N).
\end{align*}
Hence, it suffices to show only 
\begin{align*}
\sup_{n,x} E_x ((X_{\xi_n} - X_{n})_+ |\tilde {\cal F}_n, X_n > N) \le M < \infty \qquad (a.s.)
\end{align*}
In other words, it is sufficient to establish for $n=0$ that 
\begin{align*}
\sup_{x>N} E_x (X_{\xi_0} - x)_+  \le M < \infty.
\end{align*}

~

With the same notations
$e_i = 1(X_{i+1}\ge X_i)$, $\bar e_i = 1(X_{i+1}<X_i)$, $\Delta X_i = X_{i+1} - X_i$, $\ell_n^i = \bar e_i \times \prod_{k=n}^{i-1} e_k$ we have, 
\begin{align*}
E_x(X_{\xi_n} - X_{n})_+ = E_x\sum_{i=n+1}^{\infty} \ell_n^i (X_i - X_n)
= \sum_{i=n+1}^{\infty} E_x\ell_n^i (X_i - X_n)
\end{align*}
(assuming that the latter sum converges; note that all its terms are non-negative).
Further, 
\begin{align*}
E_x\ell_n^i(X_i - X_n) = E_x(\ell_n^i\sum_{j=n}^{i-1}\Delta X_j)
= \sum_{j=n}^{i-1} E_x\ell_n^i\Delta X_j
\end{align*}
For each single term in this sum we have ($n\le j\le i-1$)
\begin{align*}
E_x\ell_n^i\Delta X_j
= E_x E_{{\cal F}_{j+1}}(\prod_{k=n}^{j} e_k) \Delta X_j (\prod_{k'=j+1}^{i-1} e_{k'}) = E (\prod_{k=n}^{j} e_k) \Delta X_j E_{{\cal F}_{j+1}}(\prod_{k'=j+1}^{i-1} e_{k'})
 \\\\
 = E_x (\prod_{k=n}^{j} e_k) \Delta X_j E_{\tilde {\cal F}_{j+1}}(\prod_{k'=j+1}^{i-1} e_{k'})
\stackrel{(A3)}\le E_x (\prod_{k=n}^{j} e_k) \Delta X_j \times q^{i-j-1}
 \\\\
= q^{i-j-1} E_x (\prod_{k=n}^{j-1} e_k) E_{{\cal F}_j}e_j \Delta X_j 
= q^{i-j-1} E_x (\prod_{k=n}^{j-1} e_k) E_{\tilde {\cal F}_j}e_j \Delta X_j 
 \\\\
\stackrel{(A5)}\le  M_1 q^{i-j-1} E_x (\prod_{k=n}^{j-1} e_k)
\le M_1 q^{i-j-1} q^{j-n} = M_1 q^{i-n-1}.
\end{align*}
Hence, 
\begin{align*}
E_x\ell_n^i(X_i - X_n) = \sum_{j=n}^{i-1} E_x\ell_n^i\Delta X_j 
\le \sum_{j=n}^{i-1} M_1 q^{i-n-1} 
= (i-n)M_1 q^{i-n-1}
\end{align*}
and so
\begin{align*}
E_x(X_{\xi_n} - X_{n})_+
\le M_1 \sum_{i=n+1}^{\infty} (i-n) q^{i-n-1}=:M_4 < \infty,
\end{align*}
as required. Lemma \ref{lem3}
is proved.  \hfill QED

\fi

\section{Main results}

\begin{theorem}
\label{thm1}
Under the assumptions (A1) -- (A5) there exist constants $C_1, C_2>0$ such that 
\begin{equation}\label{tauinfty}
E_x \tau
\le x + C_1.
\end{equation}
and there exist constants $C_2, C_3>0$ such that 
\begin{equation}\label{gammainfty}
E_x \gamma \le C_2 x + C_3.
\end{equation}
Here $C_1 \le  \frac{M_4\bar q}{1-\bar q}$.

\end{theorem}

\begin{corollary}\label{cor}
The process $X_n$ has a   stationary measure.
\end{corollary}


\ifproofs

\section{Proof of theorem \ref{thm1}}
{\bf 0}. First of all let us state the idea of the proof. We will establish the property of 
recurrence towards the interval $[0,N]$ due to the recurrence assumptions, which property holds true despite the non-markovian behaviour. Further, inside $[0,N]$ coupling holds true on each step with a positive probability bounded away from zero on the jump up (or stay); after such a coupling, the process does not remember its past given the present before it started falling down. Hence, de-coupling is not possible. 

Formally, let us make the process (strong) Markov by extending its state space. For this aim it suffices to define
\begin{eqnarray}\label{markovn}
Y_n := X_n 1(X_n\ge X_{n-1}) + (X_n, \ldots X_{\zeta_n})^T 1(X_n < X_{n-1}) \equiv (X_n, \ldots X_{\zeta_n\wedge n})^T 
\end{eqnarray}
(here $T$ stands for the transposition; recall that $\zeta_n < n$ in case of $X_n < X_{n-1}$; in any case, the vector $Y_n$ is of a finite, but variable  dimension which is random). 

~

\noindent
{\bf 1}. Recurrence. Due to (\ref{prod}), from {\em any} state $y>N$ there is a positive probability to attain the set $[0,N]$ in a single monotonic fall down with no stopovers with a probability no less than $\bar\kappa_\infty$. The time required for such a monotonic trajectory from $y$ to $[0,N]$  is no more than $y-N-1$. However, other scenarios are possible with stopovers and temporary runs up. Hence, to evaluate the expected value of $\tau$ some calculus is needed.

~

Let us establish the bound (\ref{tauinfty}).
\begin{equation}\label{etau}
E_x\tau \le x + C.
\end{equation}
If $x\le N$, then $\tau=0$ and the bound is trivial. Let $x>N$. Recall that slightly abusing notations we only write down the initial position $x$, while in fact there might be some non-trivial prehistory $\tilde {\cal F}_0$. The process may start descending straight away, or after several steps up (or after staying at state $x$ for some time). In the latter case the position $X_{\xi_0}$ from which the descent starts admits the bound
$$
(E_x X_{\xi_0}-x)_+ 
\le M_4 
$$
(see lemma \ref{lem3}).

~

\underline{Case I: at $t=0$ the process is falling down.}\\
Let us define stopping times
$$
t_0=T_0 = 0, \, T_1 = \chi_{t_0}, \, t_1 = \xi_{T_1},\, T_2 =\chi_{t_1}, \, 
t_2 = \xi_{T_2},\, T_3 =\chi_{t_2}, \ldots
$$
{\em In words, $T_i$ is the end of the next after $t_{i-1}$ partial fall; $t_i$ is the end of the next after $T_i$ run up. There might be a.s. finitely many excursions down and up, and the last fall down will finish at $[0,N]$.}

Let us recall that 
$$
\xi_n : = \max(k\ge n: \text{all increments}\; \Delta X_i\ge 0, \, \forall \, n\le i\le k)\vee n, 
$$
and
$$
\chi_n : = \max(k\ge n: \text{all increments}\; \Delta X_i < 0, \, \forall \, n\le i\le k)\vee n.
$$
We have $\forall x>N$
$$
E_x (\xi_0 - 0) \le \sum_i i q^i = : M_2.
$$
and  $\forall x>N$
$$
E_x (\chi_0 - 0) \le M_3.
$$
Note that 
$$
T_i - t_{i-1} \le X_{t_{i-1}}.
$$

Denote by $A_i$ ($i\ge 1$) the event of precisely $i-1$ unsuccessful attempts to descend to the floor $[0,N]$, after which on the $i$th attempt it does attain the floor; by $B_j$ let us denote $j$th unsuccessful attempt to fall down until reaching the floor $[0,N]$ (probability that it is unsuccessful is less that $\bar q < 1$); $B_j^c$ is the event where the $j$th fall down is successful. Then we have $\tau = T_{i}$ on $A_i = (\bigcap_{1\le j\le i-1}B_j)\bigcap B^c_{i}$. The probability of $A_i$ does not exceed $\bar q^{i-1}$. (Recall, $\bar q = 1 - \bar \kappa_\infty$.) So, we estimate, 
\begin{align*}
 E_x\tau = \sum_{i\ge 1} E_x\tau 1(A_i)
 \le  \sum_{i\ge 1} E_x 1(A_i) T_{i}   
 =  \sum_{i\ge 1} E_x 1(\bigcap_{1\le j\le i-1}B_j)\bigcap B^c_{i}) T_{i}   
 \\ \\
\stackrel{2}{=}  \sum_{i\ge 1} E_x \left(\prod_{1\le j\le i-1}1(B_j)\right) 1(B^c_{i}) T_{i}   
\stackrel{3}{=}  \sum_{i\ge 1} E_x E_{{\cal F}_{t_{i-1}}}\left(\prod_{1\le j\le i-1}1(B_j)\right) 1(B^c_{i}) T_{i}   
 \\ \\
\stackrel{4}{=}  \sum_{i\ge 1} E_x \left(\prod_{1\le j\le i-1}1(B_j)\right) E_{{\cal F}_{t_{i-1}}}1(B^c_{i}) T_{i}   
 \\ \\
\stackrel{5}{=}  \sum_{i\ge 1} E_x \left(\prod_{1\le j\le i-1}1(B_j)\right) E_{{\cal F}_{t_{i-1}}}1(B^c_{i}) (T_{i} - t_{i-1} + t_{i-1}) 
 \\   \\
\stackrel{6}{=}  \sum_{i\ge 1} E_x \left(\prod_{1\le j\le i-1}1(B_j)\right) 1(B^c_{i}) t_{i-1} + \sum_{i\ge 1} E_x \left(1\prod_{1\le j\le i-1}(B_j)\right) E_{{\cal F}_{t_{i-1}}} 1(B^c_{i}) (T_{i} - t_{i-1}) 
 \\ \\
\stackrel{7}{\le}  \sum_{i\ge 1} E_x \left(\prod_{1\le j\le i-1}1(B_j)\right) 1(B^c_{i}) t_{i-1} + \sum_{i\ge 1} E_x \left(1\prod_{1\le j\le i-1}(B_j)\right) E_{{\cal F}_{t_{i-1}}} 1(B^c_{i}) X_{t_{i-1}} 
 \\ \\
\stackrel{8}{\le} 
\displaystyle \sum_{i\ge 1} E_x \left(\prod_{1\le j\le i-1}1(B_j)\right) 1(B^c_{i}) t_{i-1} + \sum_i E_x \left(\prod_{1\le j\le i-1}1(B_j)\right) 1(B^c_{i})X_{t_{i-1}}.
\end{align*}
Note that $B_j \in {\cal F}_{T_j}$.
We are going to show that 
\begin{align}\label{1des}
\sum_{i\ge 1} E_x \left(\prod_{1\le j\le i-1}1(B_j)\right)  t_{i-1} \le C
\end{align}
and 
\begin{align}\label{2des}
\sum_{i\ge 1} E_x \left(\prod_{1\le j\le i-1}1(B_j)\right) 1(B^c_{i}) X_{t_{i-1}} \le x + C.
\end{align}

~

\noindent
{\bf Step 1.}
$$
t_{i-1} = (t_{i-1} - T_{i-1}) + (T_{i-1} - t_{i-2}) + ... + (T_1 - t_0) + (t_0 - T_0).
$$
We have
\begin{align*}
E_x \left(\prod_{1\le j\le i-1}1(B_j)\right)(t_{i-1} - T_{i-1})
=E_x E_{{\cal F}_{T_{i-1}}} \left(\prod_{1\le j\le i-1}1(B_j)\right)(t_{i-1} - T_{i-1})
 \\\\
= E_x \left(\prod_{1\le j\le i-1}1(B_j)\right) E_{{\cal F}_{T_{i-1}}} (t_{i-1} - T_{i-1}) 
\stackrel{lemma\, 1}\le M_2 E_x \left(\prod_{1\le j\le i-1}1(B_j)\right) \le M_2 \bar q^{i-1};
\end{align*}
also, 
\begin{align*}
E_x \left(\prod_{1\le j\le i-1}1(B_j)\right)(T_{i-1} - t_{i-2})
=E_x E_{{\cal F}_{t_{i-2}}} \left(\prod_{1\le j\le i-1}1(B_j)\right)(T_{i-1} - t_{i-2})
 \\\\
= E_x \left(\prod_{1\le j\le i-2}1(B_j)\right) E_{{\cal F}_{t_{i-2}}}  1(B_{i-1}) (T_{i-1} - t_{i-2})
\stackrel{lemma\, 2}\le M_3 E_x \left(\prod_{1\le j\le i-2}1(B_j)\right)
\le M_3 \bar q ^{i-2};
\end{align*}
further, 
\begin{align*}
E_x \left(\prod_{1\le j\le i-1}1(B_j)\right)(t_{i-2} - T_{i-2})
=E_x  \left(\prod_{1\le j\le i-2}1(B_j)\right)(t_{i-2} - T_{i-2}) 
E_{{\cal F}_{t_{i-1}}} 1(B_{i-1})
 \\\\
\le \bar q\, E_x \left(\prod_{1\le j\le i-2}1(B_j)\right)  (t_{i-2} - T_{i-2}) 
\le \bar q M_2 \bar q^{i-2} = M_2 \bar q^{i-1},
\end{align*}
and
\begin{align*}
E_x \left(\prod_{1\le j\le i-1}1(B_j)\right)(T_{i-2} - t_{i-3})
=E_x E_{{\cal F}_{T_{i-2}}} \left(\prod_{1\le j\le i-1}1(B_j)\right)(T_{i-2} - t_{i-3})
 \\\\
= E_x \left(\prod_{1\le j\le i-2}1(B_j)\right) (T_{i-2} - t_{i-3}) E_{{\cal F}_{T_{i-2}}} 1(B_{i-1}) 
 \\\\
\le \bar q E_x \left(\prod_{1\le j\le i-2}1(B_j)\right) (T_{i-2} - t_{i-3}) 
\le \bar q M_3 \bar q ^{i-3} = M_3 \bar q^{i-2};
\end{align*}
etc. By induction we obtain
\begin{align*}
E_x  \left(\prod_{1\le j\le i-1}1(B_j)\right) t_{i-1}
\le i M_2 \bar q^{i-1} + (i-1) M_3\bar q^{i-2}.
\end{align*}
Hence, the first desired inequality (\ref{1des}) is true, 
\begin{align*}
\sum_i E_x \left(\prod_{1\le j\le i-1}1(B_j)\right)  t_{i-1} 
\le M_2 \sum_{i\ge 1} (i-1) \bar q^{i-1} +  M_3 \sum_{i\ge 2} (i-2) \bar q^{i-1} =: C < \infty.
\end{align*}

~

\noindent
{\bf Step 2.} Note that $X_{t_{j}} \ge X_{T_{j}}$, so that 
$X_{t_j} - X_{t_{j-1}} \le X_{t_j} - X_{T_{j}}$. Also, in the case under the consideration $X_{t_{0}} = x$. Hence, we have, 
\begin{align*}
X_{t_{i-1}} = (X_{t_{i-1}} - X_{t_{i-2}}) + \ldots + (X_{t_{1}} - X_{t_{0}}) + (X_{t_{0}} - x) + x.
\end{align*}
So, 
\begin{align*}
E_x \left(\prod_{1\le j\le i-1}1(B_j)\right) 1(B^c_{i}) X_{t_{i-1}} 
 \\\\
=  E_x \left(\prod_{1\le j\le i-1}1(B_j)\right) 
1(B^c_{i}) \left(x + (X_{t_{0}}-x)+ \sum_{k=1}^{i-1} (X_{t_{k}} - X_{t_{k-1}})\right)
 \\\\
\le E_x \left(\prod_{1\le j\le i-1}1(B_j)\right) 
1(B^c_{i}) \left(x + \sum_{k=1}^{i-1} (X_{t_{k}} - X_{T_{k}})\right)
 \\\\
= x E_x\prod_{1\le j\le i-1}1(B_j)1(B^c_{i})  + E_x \left(\prod_{1\le j\le i-1}1(B_j)\right) 1(B^c_{i}) 
\sum_{k=1}^{i-1} (X_{t_{k}} - X_{T_{k}})
 \\\\
\le x E_x\prod_{1\le j\le i-1}1(B_j)1(B^c_{i})  + E_x \left(\prod_{1\le j\le i-1}1(B_j)\right) 
\sum_{k=1}^{i-1} (X_{t_{k}} - X_{T_{k}}).
\end{align*}
For any $1\le k \le i-1$ we estimate
\begin{align*}
E_x \left(\prod_{1\le j\le i-1}1(B_j)\right) 
(X_{t_{k}} - X_{T_{k}}) 
= E_x E_{{\cal F}_{t_{k}}}\left(\prod_{1\le j\le i-1}1(B_j)\right) 
(X_{t_{k}} - X_{T_{k}})
 \\\\
= E_x \left(\prod_{1\le j\le k}1(B_j)\right) (X_{t_{k}} - X_{T_{k}})
E_{{\cal F}_{t_{k}}}\left(\prod_{k+1\le j\le i-1}1(B_j)\right) 
 \\\\
\le E_x \left(\prod_{1\le j\le k}1(B_j)\right) (X_{t_{k}} - X_{T_{k}})
\bar q^{i-k-1} 
= \bar q^{i-k-1} E_x \left(\prod_{1\le j\le k}1(B_j)\right) 
E_{{\cal F}_{T_{k}}}(X_{t_{k}} - X_{T_{k}})
 \\\\
\stackrel{lemma \, 3} \le M_4 \bar q^{i-k-1} E_x \left(\prod_{1\le j\le k}1(B_j)\right) \le  M_4 \bar q^{i-k-1 + k} =  M_4 \bar q^{i-1}.
\end{align*}
Therefore, since $1 = \sum_i \left(\prod_{1\le j\le i-1}1(B_j)\right) 1(B^c_{i})$ a.s., we get
\begin{align*}
\sum_i E_x \left(\prod_{1\le j\le i-1}1(B_j)\right)1(B^c_{i})  X_{t_{i-1}} 
 \\\\
\le x E_x \sum_i \prod_{1\le j\le i-1}1(B_j)1(B^c_{i})   + M_4 \sum_i i\bar q^{i-1}
\le x + \frac{M_4}{1-\bar q}.
\end{align*}
This shows (\ref{2des}), as required.



~

\underline{Case II: at $t=0$ the process is going up.} 
Let us define stopping times 
$$
T_0 = 0, \, t_0= \xi_{0}, \, T_1 = \chi_{t_0}, \, t_1 = T_1 + \xi_{T_1},\, T_2 =t_1+ \chi_{t_1}, \, 
\ldots
$$
{($T_i$ is the end of the next after $t_{i-1}$ partial fall; $t_i$ is the end of the next after $T_i$ run up. There might be a.s. finitely many excursions down and up, and the last fall down will finish at $[0,N]$.)} We have, 
$$
t_{i-1} = (t_{i-1} - T_{i-1}) + (T_{i-1} - t_{i-2}) + ... + (T_1 - t_0)
+ (t_0 - T_0).
$$
So, we  estimate
\begin{align*}
 E_x\tau = \sum_{i\ge 1} E_x\tau 1(A_i) 
 =  \sum_{i\ge 1} E_x 1(A_i) T_{i}   
 =  \sum_{i\ge 1} E_x 1(\bigcap_{1\le j\le i-1}B_j)\bigcap B^c_{i}) T_{i}   
 \\ \\
\stackrel{2}{=}  \sum_{i\ge 1} E_x \left(\prod_{1\le j\le i-1}1(B_j)\right) 1(B^c_{i}) T_{i}   
\stackrel{3}{=}  \sum_{i\ge 1} E_x E_{{\cal F}_{t_{i-1}}}\left(\prod_{1\le j\le i-1}1(B_j)\right) 1(B^c_{i}) T_{i}   
 \\ \\
\stackrel{4}{=}  \sum_{i\ge 1} E_x \left(\prod_{1\le j\le i-1}1(B_j)\right) E_{{\cal F}_{t_{i-1}}}1(B^c_{i}) T_{i}   
 \\ \\
\stackrel{5}{=}  \sum_{i\ge 1} E_x \left(\prod_{1\le j\le i-1}1(B_j)\right) E_{{\cal F}_{t_{i-1}}}1(B^c_{i}) (T_{i} - t_{i-1} + t_{i-1}) 
 \\   \\
\stackrel{6}{=}  \sum_{i\ge 1} E_x \left(\prod_{1\le j\le i-1}1(B_j)\right) 1(B^c_{i}) t_{i-1} + \sum_{i\ge 1} E_x \left(1\prod_{1\le j\le i-1}(B_j)\right) E_{{\cal F}_{t_{i-1}}} 1(B^c_{i}) (T_{i} - t_{i-1}) 
 \\ \\
\stackrel{7}{\le}  \sum_{i\ge 1} E_x \left(\prod_{1\le j\le i-1}1(B_j)\right) 1(B^c_{i}) t_{i-1} + \sum_{i\ge 1} E_x \left(1\prod_{1\le j\le i-1}(B_j)\right) E_{{\cal F}_{t_{i-1}}} 1(B^c_{i}) X_{t_{i-1}} 
 \\ \\
\stackrel{8}{\le} 
\displaystyle \sum_{i\ge 1} E_x \left(\prod_{1\le j\le i-1}1(B_j)\right) 1(B^c_{i}) t_{i-1} + \sum_{i\ge 1} E_x \left(\prod_{1\le j\le i-1}1(B_j)\right) 1(B^c_{i})X_{t_{i-1}}
\end{align*}
Note that $B_j \in {\cal F}_{T_j}$.
We are going to show that 
\begin{align}\label{1des2}
\sum_{i\ge 1} E_x \left(\prod_{1\le j\le i-1}1(B_j)\right)  t_{i-1} \le C
\end{align}
and 
\begin{align}\label{2des2}
\sum_{i\ge 1} E_x \left(\prod_{1\le j\le i-1}1(B_j)\right) 1(B^c_{i}) X_{t_{i-1}} \le x + C.
\end{align}

~

\noindent
{\bf Step 3.}
We have
\begin{align*}
E_x \left(\prod_{1\le j\le i-1}1(B_j)\right)(t_{i-1} - T_{i-1})
=E_x E_{{\cal F}_{T_{i-1}}} \left(\prod_{1\le j\le i-1}1(B_j)\right)(t_{i-1} - T_{i-1})
 \\\\
= E_x \left(\prod_{1\le j\le i-1}1(B_j)\right) E_{{\cal F}_{T_{i-1}}} (t_{i-1} - T_{i-1}) 
\stackrel{lemma\, 1}\le M_2 E_x \left(\prod_{1\le j\le i-1}1(B_j)\right) \le M_2 \bar q^{i-1},
\end{align*}
and
\begin{align*}
E_x \left(\prod_{1\le j\le i-1}1(B_j)\right)(T_{i-1} - t_{i-2})
=E_x E_{{\cal F}_{t_{i-2}}} \left(\prod_{1\le j\le i-1}1(B_j)\right)(T_{i-1} - t_{i-2})
 \\\\
= E_x \left(\prod_{1\le j\le i-2}1(B_j)\right) E_{{\cal F}_{t_{i-2}}}  1(B_{i-1}) (T_{i-1} - t_{i-2})
\stackrel{lemma\, 2}\le M_3 E_x \left(\prod_{1\le j\le i-2}1(B_j)\right)
\le M_3 \bar q ^{i-2};
\end{align*}
further, 
\begin{align*}
E_x \left(\prod_{1\le j\le i-1}1(B_j)\right)(t_{i-2} - T_{i-2})
=E_x  \left(\prod_{1\le j\le i-2}1(B_j)\right)(t_{i-2} - T_{i-2}) 
E_{{\cal F}_{t_{i-1}}} 1(B_{i-1})
 \\\\
\le \bar q\, E_x \left(\prod_{1\le j\le i-2}1(B_j)\right)  (t_{i-2} - T_{i-2}) 
\le \bar q M_2 \bar q^{i-2} = M_2 \bar q^{i-1};
\end{align*}
\begin{align*}
E_x \left(\prod_{1\le j\le i-1}1(B_j)\right)(T_{i-2} - t_{i-3})
=E_x E_{{\cal F}_{T_{i-2}}} \left(\prod_{1\le j\le i-1}1(B_j)\right)(T_{i-2} - t_{i-3})
 \\\\
= E_x \left(\prod_{1\le j\le i-2}1(B_j)\right) (T_{i-2} - t_{i-3}) E_{{\cal F}_{T_{i-2}}} 1(B_{i-1})
 \\\\ 
\le \bar q E_x \left(\prod_{1\le j\le i-2}1(B_j)\right) (T_{i-2} - t_{i-3}) 
\le \bar q M_3 \bar q ^{i-3} = M_3 \bar q^{i-2};
\end{align*}
etc. By induction we obtain
\begin{align*}
E_x  \left(\prod_{1\le j\le i-1}1(B_j)\right) t_{i-1}
\le i M_2 \bar q^{i-1} + (i-1) M_3\bar q^{i-2}.
\end{align*}
Hence, the first desired inequality (\ref{1des2}) is true, 
\begin{align*}
\sum_{i\ge 1} E_x \left(\prod_{1\le j\le i-1}1(B_j)\right)  t_{i-1} 
\le M_2 \sum_{i\ge 1} i \bar q^{i-1} +  M_3 \sum_{i\ge 2} (i-1) \bar q^{i-1} =: C < \infty.
\end{align*}

~

\noindent
{\bf Step 4.} 
Note that  $X_{T_{0}} = x$, and 
$$
X_{t_{i-1}} \le x + \sum_{j=1}^{i-1} (X_{t_{j}} - X_{T_{j}}).
$$
So, we have, 
\begin{align*}
E_x \left(\prod_{1\le j\le i-1}1(B_j)\right) 1(B^c_{i}) X_{t_{i-1}} 
\le E_x \left(\prod_{1\le j\le i-1}1(B_j)\right) 
1(B^c_{i}) \left(x + \sum_{j=1}^{i-1} (X_{t_{j}} - X_{T_{j}})\right)
 \\\\
= x E_x\prod_{1\le j\le i-1}1(B_j)1(B^c_{i})  + E_x \left(\prod_{1\le j\le i-1}1(B_j)\right) 1(B^c_{i}) 
\sum_{j=1}^{i-1} (X_{t_{j}} - X_{T_{j}})
 \\\\
\le x E_x\prod_{1\le j\le i-1}1(B_j)1(B^c_{i})  + E_x \left(\prod_{1\le j\le i-1}1(B_j)\right) 
\sum_{k=1}^{i-1} (X_{t_{k}} - X_{T_{k}}).
\end{align*}
For any $1\le k \le i-1$ we estimate
\begin{align*}
E_x \left(\prod_{1\le j\le i-1}1(B_j)\right) 
(X_{t_{k}} - X_{T_{k}}) 
= E_x E_{{\cal F}_{t_{k}}}\left(\prod_{1\le j\le i-1}1(B_j)\right) 
(X_{t_{k}} - X_{T_{k}})
 \\\\
= E_x \left(\prod_{1\le j\le k}1(B_j)\right) (X_{t_{k}} - X_{T_{k}})
E_{{\cal F}_{t_{k}}}\left(\prod_{k+1\le j\le i-1}1(B_j)\right) 
 \\\\
\le E_x \left(\prod_{1\le j\le k}1(B_j)\right) (X_{t_{k}} - X_{T_{k}})
\bar q^{i-k-1} 
= \bar q^{i-k-1} E_x \left(\prod_{1\le j\le k}1(B_j)\right) 
E_{{\cal F}_{T_{k}}}(X_{t_{k}} - X_{T_{k}})
 \\\\
\stackrel{lemma \, 3} \le M_4 \bar q^{i-k-1} E_x \left(\prod_{1\le j\le k}1(B_j)\right) \le  M_4 \bar q^{i-k-1 + k} =  M_4 \bar q^{i-1}.
\end{align*}
Therefore, since $1 = \sum_{i\ge 1} \left(\prod_{1\le j\le i-1}1(B_j)\right) 1(B^c_{i})$ a.s., we get
\begin{align*}
\sum_{i\ge 1} E_x \left(\prod_{1\le j\le i-1}1(B_j)\right)1(B^c_{i})  X_{t_{i-1}} 
 \\\\
\le x E_x \sum_{i\ge 1} \prod_{1\le j\le i-1}1(B_j)1(B^c_{i})   + M_4 \sum_i\bar q^{i-1}
\le x + \frac{M_4}{1-\bar q}.
\end{align*}
This shows (\ref{2des2}), as required. In both cases I and II 
the bound (\ref{tauinfty}) is proved.

~

\noindent
{\bf Step 5.} Let us establish the bound (\ref{gammainfty}). 
Recall the notations introduced earlier after the assumptions:
$$
\tau=\tau^1 
:= \inf(t\ge 0: X_t\le N); \quad \gamma = \gamma^1:=\inf(t\ge \tau: X_{t-1}\le X_t = N),
$$
and
\begin{align*}
 T^{n}:= \inf(t> \tau^n: X_t > N), \quad 
 \tau^{n+1}:= \inf(t > T^n: X_t\le N), \quad n\ge 1,
\end{align*}
and $T^0:=0$. Also, let
\begin{align*}
  \gamma^{n+1}:= \inf(t > \gamma^n: X_{t-1}\le X_t = N).
\end{align*}

We have due to the assumption (A5)
$$
E_x X_{T^n}\le C, \quad n\ge 1; \quad \text{also,} \quad E_x X_{T^0}=x.
$$
Therefore, by virtue of the bound (\ref{tauinfty}) we have,
$$
E_x(\tau^{n+1}-T^n) = E_x E_x(\tau^{n+1}-T^n |\tilde{\cal F}_{T^n})
\le E_x X_{T^n} + C \le C, \quad n\ge 1,
$$
and
$$
E_x(\tau^{1}-T^0) \le x+C, \quad n = 0.
$$
Also, due to the assumptions there exists $p\in (0,1)$ such that 
$$
P_x(\gamma >T^n)\le p^n \quad \Longleftrightarrow \quad 
P_x(\gamma \le T^n)\ge 1-p^n, \quad n\ge 1.
$$

Also, 
\begin{align*}
E_x(T^n - \tau^{n}) \le C, \quad n \ge 1.
\end{align*}

Thus, also 
\begin{align*}
E(T^{n+1} - T^{n}) = E(T^{n+1} - \tau^{n+1} + \tau^{n+1}- T^{n})  \le C.
\end{align*}
Moreover, 
\begin{align*}
E(T^{n+1} - T^{n} |{\cal F}^X_{T^n})   \le C.
\end{align*}
It follows by induction that 
\begin{align*}
E T^{n} \le C n + x.
\end{align*}

So, we estimate 
\begin{align*}
 E_x \gamma = \sum_{n\ge 0}^{} E_x\gamma 1(T^n < \gamma \le T^{n+1}) \le \sum_{n\ge 0}^{} E_xT^{n+1} 1(T^n < \gamma \le T^{n+1})
  \\\\
= \sum_{n\ge 0}^{} E_x E_x(T^{n+1} 1(T^n < \gamma \le T^{n+1}) |{\cal F}^X_{T^n}) 
\\\\
= \sum_{n\ge 0}^{} E_x E_x((T^{n}+T^{n+1}-T^{n}) 1(T^n < \gamma \le T^{n+1}) |{\cal F}^X_{T^n}) 
\\\\
\le \sum_{n\ge 0}^{} E_x E_x((T^{n}+T^{n+1}-T^{n}) 1(T^n < \gamma) |{\cal F}^X_{T^n}) 
\\\\
= \sum_{n\ge 0}^{} E_x T^n 1(T^n < \gamma) 
+ \sum_{n\ge 0}^{} E_x 1(T^n < \gamma)  \underbrace{E_x((T^{n+1}-T^{n}) |{\cal F}^X_{T^n})}_{\le C + x 1(n=0)}.
\end{align*}
Further, with any integer $M>0$, denoting $E_x T^n 1(T^n < \gamma)=: d_n$, we have (note that $d_0=0$), 
\begin{align*}
\sum_{n = 0}^{M} \underbrace{E_x T^n 1(T^n < \gamma)}_{=: d_n}
= \sum_{n = 0}^{M}  E_x (T^{n-1} + T^{n} - T^{n-1}) 1(T^n < \gamma) 1(T^{n-1} < \gamma) 
 \\\\
= d_0+\sum_{n = 1}^{M}  E_x T^{n-1}1(T^n < \gamma) 1(T^{n-1} < \gamma) + \sum_{n = 1}^{M}  E_x (T^{n} - T^{n-1}) 1(T^n < \gamma) 1(T^{n-1} < \gamma) 
 \\\\
= d_0+\sum_{n = 1}^{M}  E_x T^{n-1} 1(T^{n-1} < \gamma)  \underbrace{E_x(1(T^n < \gamma) | {\cal F}^X_{T^{n-1}})}_{\le p}
 \\\\
+ \sum_{n = 1}^{M}  E_x 1(T^{n-1} < \gamma)  E_x ((T^{n} - T^{n-1}) 1(T^n < \gamma) |  {\cal F}^X_{T^{n-1}})
  \\\\
\le d_0+\sum_{n = 1}^{M}  p d_{n-1}
+ \sum_{n = 1}^{M}  E_x 1(T^{n-1} < \gamma)  E_x ((T^{n} - \tau^n + \tau^n - T^{n-1}) 1(T^n < \gamma) |  {\cal F}^X_{T^{n-1}})
  \\\\
= \sum_{n = 1}^{M}  p d_{n-1}
+ \sum_{n = 1}^{M}  E_x 1(T^{n-1} < \gamma)  [E_x ((T^{n} - \tau^n) 1(T^n < \gamma) |  {\cal F}^X_{T^{n-1}}) 
 \\\\
+  \underbrace{E_x ((\tau^n - T^{n-1}) 1(T^n < \gamma) |  {\cal F}^X_{T^{n-1}}))}_{\le Cp + x1(n=1)}].
\end{align*}
We have, 
\begin{align*}
\sum_{n = 1}^{M}  E_x 1(T^{n-1} < \gamma)   \underbrace{E_x ((\tau^n - T^{n-1}) 1(T^n < \gamma) |  {\cal F}^X_{T^{n-1}}))}_{\le Cp + x1(n=1)}
 \\\\
 \le C+x+Cp \sum_{n = 1}^{M} E_x 1(T^{n-1} < \gamma)  
 \le C+ x+C \sum_{n = 0}^{M-2} p^{n}\le C+x. 
\end{align*}
Further, 
\begin{align*}
 \sum_{n = 1}^{M}  E_x 1(T^{n-1} < \gamma)  [E_x ((T^{n} - \tau^n) 1(T^n < \gamma) |  {\cal F}^X_{T^{n-1}})] 
  \\\\
=  \sum_{n = 1}^{M}  E_x 1(T^{n-1} < \gamma) E_x [ E_x\{(T^{n} - \tau^n) 1(T^n < \gamma)| {\cal F}^X_{\tau^{n}}\} |  {\cal F}^X_{T^{n-1}}] 
 \\\\
\le \sum_{n = 1}^{M}  E_x 1(T^{n-1} < \gamma) E_x [\underbrace{E_x\{(T^{n} - \tau^n) | {\cal F}^X_{\tau^{n}}\}}_{\le C} |  {\cal F}^X_{T^{n-1}}] 
 \\\\
\le C  \sum_{n = 1}^{M}  E_x 1(T^{n-1} < \gamma) 
\le C \sum_{n\ge 0}^{} p^{n-1} \le C.
\end{align*}
Thus,
\begin{align*}
\sum_{n = 0}^{M}  d_n
\le p (C+\sum_{n = 0}^{M-1}  d_n) + C + x, 
\end{align*}
which  implies by the monotone convergence theorem that 
\begin{align*}
\sum_{n = 0}^{\infty} d_n\le  C(1 + x)
\end{align*}
and
\begin{align*}
E_x\gamma \le \sum_{n\ge 0}^{} d_n + C \le  C+Cx. 
\end{align*}
The bound (\ref{gammainfty}) is justified and the proof of the theorem is completed. \hfill 
\text{QED}




\section{Proof of Corollary \ref{cor}}
The existence of an invariant measure for the process $Y$ follows from the Harris -- Khasminskii principle via the formula 

\begin{align*}
\mu^Y (A) := c \, E_{N-} \sum_{n=1}^{\gamma} 1(Y_n \in A),
\end{align*}
where $c$ is the normalising constant, $A$ is any measurable set in the state space of the process $Y$, and by $E_{N-}$ we understand the initial condition $X_0=N$ with any preceding fictitious state $X_{-1}\le N$. By the assumptions, the  distribution of $X_1$ only depends on $X_0$ given this condition. So, this state -- with the  convention of the preceding state in $[0,N]$ -- is, indeed, a regeneration point. 

To apply it to the process $X$ let us take any bounded measurable function $f(y)$ ($y = (y^1, \ldots)$),
which only depends on the first variable $y^1=x$:
\begin{align*}
\int f(y)\mu^Y (dy) := c  E_{N-} \sum_{n=1}^{\gamma} f(Y_n).
\end{align*}
The latter expression in the right hand side determines an invariant measure for the process $X$ with a notation $g(y^1):=f(y)$: 
\begin{align*}
\int g(y^1)\mu^Y (dy) := c  E_{N-} \sum_{1}^{\gamma} g(Y^1_n),
\end{align*}
where $X_n=Y^1_n$. So, the invariant measure for $X$ reads
\begin{align*}
\mu^X (A^1) := c  E_{N-} \sum_{n=1}^{\gamma} 1(Y^1_n \in A^1).
\end{align*}
The corollary is proved. \hfill QED

~

\fi

\section{Acknowledgements}
For the first author 
the part consisting of theorem \ref{thm1} and lemma \ref{lem3}  
was supported by Russian Foundation for Basic Research 
grant 20-01-00575$\_$a.


\end{document}